\documentclass[12pt,oneside]{amsproc}
\usepackage[T2A]{fontenc}
\usepackage[utf8]{inputenc}
\usepackage[russian]{babel}
\usepackage{amssymb,amsmath,amsthm}

\numberwithin{equation}{subsection}

\begin{document}
\author{А.\,А.~Владимиров, Е.\,С.~Карулина}
\title[]{О нижней априорной оценке минимального собственного значения одной задачи
Штурма--Лиувилля с граничными условиями второго типа}

\begin{abstract}
Устанавливается достижимость точной нижней грани $m_\gamma$ минимальных
собственных значений граничных задач
\begin{gather*}
	-y''+qy=\lambda y,\\
	y'(0)=y'(1)=0
\end{gather*}
при пробегании неотрицательным потенциалом $q\in L_1[0,1]$ единичной сферы
пространства $L_\gamma[0,1]$, где $\gamma\in (0,1)$. Также устанавливается
справедливость равенства $m_\gamma=1$ при $\gamma\leqslant 1-2\pi^{-2}$,
и неравенства $m_\gamma<1$ иначе.
\end{abstract}

\maketitle
\footnotetext[1]{Работа первого автора поддержана РФФИ, проект \No~13-01-00705.
Работа обоих авторов поддержана РФФИ, проект \No~14-01-31423, и РНФ,
проект \No~14-11-00754.}

\subsection{Введение}\label{par:1}
Рассмотрим задачу Штурма--Лиувилля
\begin{gather}\label{eq:eq1}
	-y''+qy=\lambda y,\\ \label{eq:eq2}
	y'(0)=y'(1)=0,
\end{gather}
где \(y\in W_1^2[0,1]\), равенство \eqref{eq:eq1} понимается в смысле
обобщенного дифференцирования, а потенциал $q$ пробегает множество
$$
	A_\gamma\rightleftharpoons\left\{q\in L_1[0,1]\;:\;
		q(x)\stackrel{\text{п.\,в.}}{\geqslant}0,\quad
		\int_0^1 q^\gamma\,dx=1\right\}.
$$
Целью настоящей статьи является уточнение некоторых свойств точной нижней грани
$m_{\gamma}\rightleftharpoons\inf_{q\in A_{\gamma}}\lambda_{1}(q)$ собственных
значений задачи \eqref{eq:eq1},~\eqref{eq:eq2} в случае $\gamma\in (0,1)$.
А именно, нами будут установлены следующие два предложения:

\medskip\subsubsection{}\label{theorem1}
{\itshape При всяком $\gamma\in (0,1)$ существует потенциал $q_*\in A_{\gamma}$,
удовлетворяющий равенству $\lambda_1(q_*)= m_{\gamma}$.
}

\medskip\subsubsection{}\label{theorem2}
{\itshape При всяком $\gamma\in (0,1-2\pi^{-2}]$ выполняется равенство $m_\gamma=1$,
а при всяком $\gamma\in (1-2\pi^{-2},1)$ выполняется неравенство $m_\gamma<1$.
}

\medskip
Все рассматриваемые в тексте функциональные пространства предполагаются
вещественными.

Отметим ряд известных в литературе результатов для оптимизационных проблем,
подобных рассматриваемой нами. В работе \cite{EK:1996} было исследовано при всех
$\gamma\in\mathbb R\setminus\{0\}$ поведение величины $m_\gamma$, отвечающей
первой граничной задаче для уравнения
\begin{equation}\label{eq:eq1'}
	-y''=\lambda qy.
\end{equation}
Та же проблема в случае второй и третьей граничных задач для уравнения
\eqref{eq:eq1'} изучена в работе \cite{Mur:2005}. Случай первой граничной задачи
для уравнения \eqref{eq:eq1} рассмотрен в работе \cite{Ez:2005}. В тех же работах
для рассматриваемых в них граничных задач изучено поведение не рассматриваемой
нами величины $M_\gamma\rightleftharpoons\sup_{q\in A_{\gamma}}\lambda_{1}(q)$.
Поведение последней при $\gamma\geqslant 1$ в первой граничной задаче для уравнения
\eqref{eq:eq1} рассматривалось также в работе \cite{VS:2003}, а при $\gamma=1$~---
в работе \cite{EK:1984}.

В случае $\gamma\geqslant 1$ рассматриваемая в настоящей статье задача
изучалась, в частности, в работах \cite{Kar:Uniti}, \cite{K_V:Tatra}.


\subsection{Переформулировка задачи}\label{par:2}
Далее мы всегда будем исходить из факта равносильности (см., например,
\cite{K_V:Tatra}) граничной задачи \eqref{eq:eq1}, \eqref{eq:eq2} спектральной
задаче для операторного пучка $T_q\colon\mathbb R\to\mathcal B(W_2^1[0,1],
W_2^{-1}[0,1])$ вида
$$
	\langle T_q(\lambda)y,y\rangle\equiv\int_0^1\bigl [|y'|^2+
		(q-\lambda)y^2\bigr]\,dx.
$$

Введем в рассмотрение параметризованное значениями $\zeta>0$ семейство множеств
$$
	M_\zeta\rightleftharpoons\{y\in W_2^1[0,1]\;:\; (\forall x\in [0,1])\quad
		y(x)\geqslant\zeta\}.
$$
При этом множество $M$ всех равномерно положительных функций класса $W_2^1[0,1]$
очевидным образом допускает представление
$$
	M=\bigcup_{\zeta>0} M_\zeta.
$$
Введем также в рассмотрение функционалы $J_\gamma\colon M\to\mathbb R$
и $G_\gamma\colon M\to\mathbb R$ вида
\begin{gather*}
	J_{\gamma}(y)\rightleftharpoons\int_0^1 (y')^2\,dx+
		\left(\int_0^1 y^{2\gamma/(\gamma-1)}\,dx
		\right)^{(\gamma-1)/\gamma},\\
	G_\gamma(y)\rightleftharpoons J_\gamma(y)\cdot\|y\|_{L_2[0,1]}^{-2}.
\end{gather*}
Имеет место следующий факт:

\medskip\subsubsection{}
{\itshape Справедливо равенство $m_\gamma=\inf_{y\in M} G_\gamma(y)$.}
\begin{proof}
Согласно неравенству Гельдера, для любых $q\in A_\gamma$ и $y\in M$ справедливы
соотношения
$$
	1=\left(\int_0^1q^\gamma\,dx\right)^{1/\gamma}\leqslant
		\left[\left(\int_0^1qy^2\,dx\right)^\gamma\cdot
		\left(\int_0^1 y^{\frac{2\gamma}{\gamma-1}}\,dx
		\right)^{1-\gamma} \right]^{1/\gamma},
$$
влекущие независимо от выбора $\lambda\in\mathbb R$ оценку $\langle T_q(\lambda)y,
y\rangle \geqslant J_\gamma(y)-\lambda\,\|y\|_{L_2[0,1]}^2$. Соответственно,
зафиксировав (в согласии с теорией Штурма) в качестве $y\in M$ собственную функцию
пучка $T_q$, отвечающую его наименьшему собственному значению $\lambda_1(q)$,
устанавливаем справедливость соотношений
$$
	G_\gamma(y)\leqslant\langle T_q(\lambda_1(q))y,y\rangle\cdot
		\|y\|_{L_2[0,1]}^{-2}+\lambda_1(q)=\lambda_1(q).
$$
Ввиду произвольности выбора потенциала $q\in A_\gamma$ это означает справедливость
оценки $m_\gamma\geqslant\inf_{y\in M}G_\gamma(y)$.

Далее, при любом выборе функции $y\in M$ для потенциала
$$
	q_*\rightleftharpoons\left(\int_0^1 y^{2\gamma/(\gamma-1)}\,dx
		\right)^{-1/\gamma}\cdot y^{2/(\gamma-1)}\in A_\gamma
$$
справедливы равенства $\langle T_{q_*}(G_\gamma(y))y,y\rangle=J_\gamma(y)-
G_\gamma(y)\,\|y\|_{L_2[0,1]}^2=0$, а тогда и оценка $\lambda_1(q_*)\leqslant
G_\gamma(y)$. Произвольность выбора функции $y\in M$ означает потому
справедливость оценки $m_\gamma\leqslant\inf_{y\in M}G_\gamma(y)$.
\end{proof}


\subsection{Доказательство предложения \ref{theorem1}}\label{par:3}
Имеет место следующий факт:

\medskip\subsubsection{}\label{prop:3.1.1}
{\itshape Пусть $\varepsilon\in (0,1)$, и пусть $y\in M$ удовлетворяет неравенству
\begin{equation}\label{neq:J_not_large}
	G_\gamma(y)<(\pi^2/4)\cdot(1-\varepsilon)^2.
\end{equation}
Тогда независимо от выбора точки $x\in[0,1]$ справедлива оценка
\begin{equation}\label{neq:y_not_small}
	\belowdisplayskip 0pt
	y(x)>\varepsilon\,\|y\|_{L_2[0,1]}.
\end{equation}
}
\begin{proof}
Предположим, что некоторая точка $a\in [0,1]$ удовлетворяет неравенству
$y(a)\leqslant\varepsilon\,\|y\|_{L_2[0,1]}$. Рассмотрим принадлежащие пространству
$W_2^1[0,1]$ функции
$$
	y_-(x)\rightleftharpoons\begin{cases}
			y(x)-y(a)&\text{при }x\leqslant a,\\
			0&\text{при }x\geqslant a,\end{cases}\qquad
	y_+(x)\rightleftharpoons\begin{cases}
			y(x)-y(a)&\text{при }x\geqslant a,\\
			0&\text{при }x\leqslant a.\end{cases}
$$
При этом очевидным образом справедливы оценки
\begin{align}
	\|y_-\|_{L_2[0,1]}^2+\|y_+\|_{L_2[0,1]}^2
		&=\|y-y(a)\|_{L_2[0,1]}^2\notag \\ \label{eq:2}
		&\geqslant(1-\varepsilon)^2\,\|y\|_{L_2[0,1]}^2.
\end{align}
С другой стороны, справедливы также соотношения
\begin{align*}
	J_\gamma(y)&\geqslant\int_0^1 (y_-')^2\,dx+\int_0^1 (y_+')^2\,dx\\
		&\geqslant(\pi^2/4)\,\|y_-\|_{L_2[0,1]}^2+(\pi^2/4)\,
		\|y_+\|_{L_2[0,1]}^2,
\end{align*}
объединяя которые с \eqref{eq:2}, получаем оценку $G_\gamma(y)\geqslant
(\pi^2/4)\cdot(1-\varepsilon)^2$.
\end{proof}

Обозначим символом $Y$ множество функций $y\in M$, для которых неравенство
\eqref{neq:y_not_small} нарушается хотя бы в одной точке отрезка $[0,1]$.
Введем также обозначение $\hat y$ для функции вида $\hat y(x)\equiv 1$.
Предложение~\ref{prop:3.1.1} гарантирует в случае $\varepsilon<1/3$
справедливость оценок
$$
	m_\gamma\leqslant G_\gamma(\hat y)\\ =1<(\pi^2/4)\cdot(1-\varepsilon)^2
		\leqslant\inf_{y\in Y} G_\gamma(y),
$$
а потому и соотношения $m_\gamma=\inf_{y\in M\setminus Y} G_\gamma(y)$. Незначительной
модификацией этого факта является следующее предложение, в формулировке которого
мы используем символ $S$ для обозначения единичной сферы пространства
$W_2^1[0,1]$:

\medskip\subsubsection{}\label{prop:3.1.2}
{\itshape Существует значение $\zeta>0$ со свойством $m_\gamma=\inf_{y\in M_\zeta
\cap S}G_\gamma(y)$.
}
\begin{proof}
Достаточно заметить, что всякая функция $y\in M\setminus Y$ со свойством
$G_\gamma(y)<3$ удовлетворяет оценкам
$$
	\|y\|_{W_2^1[0,1]}^2\leqslant [G_\gamma(y)+1]\cdot\|y\|_{L_2[0,1]}^2
		<4\,\|y\|_{L_2[0,1]}^2,
$$
гарантирующим справедливость оценок $y(x)>(\varepsilon/2)\cdot\|y\|_{W_2^1[0,1]}$.
\end{proof}

Заметим теперь, что полная непрерывность вложения пространства $W_2^1[0,1]$
в пространство $L_2[0,1]$ гарантирует существование последовательности
$\{y_n\}_{n=0}^\infty$ функций $y_n\in M_\zeta\cap S$, обладающей
свойством $\lim_{n\to\infty}G_\gamma(y_n)=m_\gamma$ и сходящейся в пространстве
$L_2[0,1]$ к некоторой функции $y_*\in L_2[0,1]$. Рассмотрим соответствующую
двойную последовательность $\{y_{n,m}\}_{(n,m)\in\mathbb N^2}$ функций вида
$$
	y_{n,m}\rightleftharpoons\dfrac{y_n+y_m}2\in M_\zeta.
$$
С учетом неравенства Гельдера и теоремы Лагранжа о конечном приращении,
указанные функции подчиняются оценкам
$$
	\left|\int_0^1 y_{n,m}^{2\gamma/(\gamma-1)}\,dx-
		\int_0^1 y_*^{2\gamma/(\gamma-1)}\,dx\right|\leqslant
		\dfrac{2\gamma\cdot\zeta^{(\gamma+1)/(\gamma-1)}}{1-\gamma}\cdot
		\|y_{n,m}-y_*\|_{L_2[0,1]}.
$$
Соответственно, справедливы также соотношения
\begin{align*}
	\liminf_{n,m\to\infty}G_\gamma(y_{n,m})&=\dfrac{\displaystyle
		\liminf_{n,m\to\infty}\|y_{n,m}\|_{W_2^1[0,1]}^2-\|y_*\|_{L_2[0,1]}^2
		+\left(\int_0^1 y_*^{2\gamma/(\gamma-1)}\,dx\right)^{(\gamma-1)/
		\gamma}}{\|y_*\|_{L_2[0,1]}^2}\\
		&\leqslant\dfrac{\displaystyle 1-\|y_*\|_{L_2[0,1]}^2
		+\left(\int_0^1 y_*^{2\gamma/(\gamma-1)}\,dx\right)^{(\gamma-1)/
		\gamma}}{\|y_*\|_{L_2[0,1]}^2}\\
		&=\lim_{n\to\infty}G_\gamma(y_n)\\ &=m_\gamma,
\end{align*}
вместе с тривиальной оценкой $\liminf_{n,m\to\infty}G_\gamma(y_{n,m})\geqslant
m_\gamma$ означающие выполнение равенства $\lim_{n,m\to\infty}\|y_{n,m}\|_{W_2^1[0,1]}
=1$. Из факта равномерной выпуклости гильбертова пространства $W_2^1[0,1]$ теперь
немедленно вытекает, что последовательность $\{y_n\}_{n=0}^\infty$ фундаментальна
в пространстве $W_2^1[0,1]$. Тем самым, предельная функция $y_*\in L_2[0,1]$
в действительности принадлежит классу $M_\zeta\cap S$ и удовлетворяет
равенству $G_\gamma(y_*)=m_\gamma$.

Итак, нами установлено существование функции $y_*\in M$, на которой достигается
минимальное по классу $M$ значение функционала $G_\gamma$. Это автоматически
означает обращение в нуль дифференциала $DG_\gamma(y_*)$, что, в свою очередь,
равносильно выполнению равенства $T_{q_*}(m_\gamma)y_*=0$, где положено
$$
	q_*\rightleftharpoons\left(\int_0^1 y_*^{2\gamma/(\gamma-1)}\,dx
		\right)^{-1/\gamma}\cdot y_*^{2/(\gamma-1)}\in A_\gamma.
$$
Доказательство предложения \ref{theorem1} тем самым завершено.


\subsection{Доказательство предложения \ref{theorem2}}\label{par:4}
Начнем с установления справедливости второй части рассматриваемого предложения:

\medskip\subsubsection{}\label{prop:4.1.1}
{\it При любом $\gamma\in (1-2\pi^{-2},1)$ выполняется неравенство $m_\gamma<1$.
}
\begin{proof}
Достаточно установить, что функция $\hat y\in M$ вида $\hat y(x)\equiv 1$ не является
минимумом функционала $G_\gamma$. Это заведомо имеет место, когда квадратичная форма
второго дифференциала функционала $G_\gamma$ в точке $\hat y$ не является
неотрицательно определенной. Для доказательства последнего, ввиду легко проверяемых
тождеств
\begin{gather*}
    DG_\gamma(y;\,v)\equiv\dfrac{2}{\|y\|_{L_2[0,1]}^2}\cdot\biggl[
		\int_0^1 y'v'\,dx+{}\kern 4cm\\
		\kern 2cm{}+\left(\int_0^1 y^{2\gamma/(\gamma-1)}\,dx
			\right)^{-1/\gamma}\cdot\int_0^1 y^{(\gamma+1)/(\gamma-1)}v
			\,dx-G_\gamma(y)\cdot\int_0^1 yv\,dx\biggr],
\end{gather*}
$$
	D^2G_\gamma(\hat y;\,u,v)\equiv\int_0^1 2u'v'\,dx+\dfrac{4}{1-\gamma}\,
			\left(\int_0^1 u\,dx\cdot\int_0^1 v\,dx-\int_0^1 uv\,dx
			\right),
$$
достаточно установить отрицательность наименьшего собственного значения
граничной задачи
\begin{gather*}
	-y''+\dfrac{2}{1-\gamma}\,\left(\int_0^1 y\,dx-y\right)=\lambda y,\\
 	y'(0)=y'(1)=0.
\end{gather*}
Однако это собственное значение есть $\pi^2-2\,(1-\gamma)^{-1}<0$.
\end{proof}

При доказательстве первой части предложения \ref{theorem2} мы будем опираться
на следующий факт:

\medskip\subsubsection{}\label{prop:notincreasing}
{\it Функция $\gamma\mapsto m_\gamma$ невозрастает на множестве $\mathbb R
\setminus\{0\}$.
}
\begin{proof}
Зафиксируем произвольную равномерно положительную функцию $q\in C[0,1]$. С учетом
неравенства Гельдера, при любом выборе значений $\gamma>0$ и $\beta>1$
справедливы оценки
$$
	\left(\int_0^1 q^\gamma\,dx\right)^{1/\gamma}=\left[\left(
		\int_0^1 q^{\beta\cdot(\gamma/\beta)}\,dx\right)^{1/\beta}
		\right]^{\beta/\gamma}\geqslant
		\left(\int_0^1 q^{\gamma/\beta}\,dx\right)^{\beta/\gamma},
$$
а при любом выборе значений $\gamma<0$ и $\beta>1$ справедливы оценки
$$
	\left(\int_0^1 q^\gamma\,dx\right)^{1/\gamma}=\left[\left(
		\int_0^1 q^{\beta\cdot(\gamma/\beta)}\,dx\right)^{1/\beta}
		\right]^{\beta/\gamma}\leqslant
		\left(\int_0^1 q^{\gamma/\beta}\,dx\right)^{\beta/\gamma}.
$$
Вместе с легко проверяемым асимптотическим соотношением
$$
	\lim_{\gamma\to 0}\left(\int_0^1 q^\gamma\,dx\right)^{1/\gamma}=
		\exp\left(\int_0^1\ln q\,dx\right)
$$
это означает неубывание функции $\gamma\mapsto\left(\int_0^1 q^{\gamma}\,
dx\right)^{1/\gamma}$ на множестве $\mathbb R\setminus\{0\}$. Соответственно,
независимо от выбора значений $\gamma,\gamma_1\in\mathbb R\setminus\{0\}$
со свойством $\gamma_1>\gamma$ для любой функции $q\in A_\gamma\cap C[0,1]$
найдется величина $C\in (0,1]$ со свойством $Cq\in A_{\gamma_1}$. Учет тривиальных
оценок $\lambda_1(Cq)\leqslant\lambda_1(q)$ и факта плотности подмножества $A_\gamma
\cap C[0,1]$ в множестве $A_\gamma\subset L_1[0,1]$ завершает доказательство.
\end{proof}

Предложение \ref{prop:notincreasing} вместе с выполняющимся независимо от выбора
$\gamma\in (0,1)$ для функции $\hat y\in M$ вида $\hat y(x)\equiv 1$ равенством
$G_\gamma(\hat y)=1$ показывают, что справедливость предложения \ref{theorem2}
достаточно проверить в случае $\gamma=1-2\pi^{-2}$. Выполнение последнего равенства
и будет предполагаться в оставшейся части статьи.

Введем в рассмотрение параметризованное величиной $\alpha>0$ семейство функций
\begin{equation} \label{eq:fa}
	f_\alpha(t)\rightleftharpoons\alpha t^{2\gamma}-t^2-1,\qquad
		t\in (0,+\infty).
\end{equation}
Имеет место следующий факт:

\medskip\subsubsection{}\label{prop:f_alpha_positive}
{\itshape Множество точек $t\in\mathbb R$ со свойством $f_\alpha(t)>0$ непусто
в точности при выполнении условия
\begin{equation} \label{eq:duonaks}
	\alpha>\gamma^{-\gamma}\,(1-\gamma)^{\gamma-1}
\end{equation}
и представляет собой в этом случае некоторый интервал $(\omega^-_\alpha,
\omega^+_\alpha)$.
}
\begin{proof}
При любом значении параметра $\alpha>0$ справедливо тождество
\begin{equation}
        f_\alpha'(t)\equiv 2\alpha\gamma t^{2\gamma-1}-2t,\label{eq:fa1}
\end{equation}
означающее, что функция $f_\alpha$ строго возрастает слева от точки $\tau_\alpha
\rightleftharpoons(\alpha\gamma)^{\frac{1}{2-2\gamma}}$ и строго убывает справа
от нее. Соответственно, область положительности функции $f_\alpha$ непуста в точности
при выполнении условия $f_\alpha(\tau_\alpha)>0$, равносильного условию
\eqref{eq:duonaks}, и представляет собой в этом случае содержащий точку $\tau_\alpha$
интервал.
\end{proof}

Кроме тождества \eqref{eq:fa1}, нами далее будут использоваться также тождества
\begin{align}
	f_\alpha''(t)&\equiv 2\alpha\gamma\,(2\gamma-1)t^{2\gamma-2}-2,\label{eq:fa2}\\
	f_\alpha'''(t)&\equiv 2\alpha\gamma(2\gamma-1)(2\gamma-2)\,
			t^{2\gamma-3}.\label{eq:fa3}
\end{align}

Предложение \ref{theorem1} показывает, что граничная задача
\begin{gather}\label{eq:1}
	-y''+y^{(\gamma+1)/(\gamma-1)}=\mu y,\\ \notag 
	y'(0)=y'(1)=0,\\ \notag\int_0^1 y^{2\gamma/(\gamma-1)}\,dx=1
\end{gather}
не может оказаться неразрешимой в случае $\mu=m_\gamma$. Непосредственно проверяется
также, что в случае $\mu<1$ соответствующее решение не может быть постоянным.
Тем самым, для завершения доказательства предложения \ref{theorem2} достаточно
установить, что период непостоянных решений уравнения \eqref{eq:1} не может быть
равен~$2$ при выполнении неравенства $\mu<1$. Поскольку, как показывают стандартные
вычисления, этот период представляет собой величину $2\mu^{-1/2}\,(1-\gamma)\,
I_0(\hat\alpha)$, где положено
\begin{gather*}
	I_\varepsilon(\alpha)\rightleftharpoons\int_{\omega^-_\alpha+
		\varepsilon}^{\omega^+_\alpha-\varepsilon}
		\dfrac{dt}{\sqrt{f_\alpha(t)}},\qquad\varepsilon\geqslant 0,\\
	\hat\alpha\rightleftharpoons\left(\dfrac{\gamma}{1-\gamma}\right)^{1-\gamma}
		\cdot\mu^{-\gamma}\cdot\left[(y')^2+\mu y^2+
		\dfrac{1-\gamma}{\gamma}\,y^{2\gamma/(\gamma-1)}\right],
\end{gather*}
то поставленная цель будет достигнута установлением справедливости следующего
предложения:

\medskip\subsubsection{}\label{prop:pi2shao}
{\itshape Функция $I_0$ не может принимать значения, меньшие величины $(1-\gamma)^{-1}$.
}

\begin{proof}
Введем обозначения
\begin{align}\label{eq:A}
	A_\alpha(t)&\rightleftharpoons t^{2\gamma}f'_\alpha(t)-
			4\gamma t^{2\gamma-1}f_\alpha(t),\\
	\psi_\alpha&\rightleftharpoons (f_\alpha')^2-
			2f_\alpha f_\alpha''.\label{eq:psi}
\end{align}
Из тождества \eqref{eq:fa3} следует, что при любой фиксации значения $\alpha_0>
\gamma^{-\gamma}\,(1-\gamma)^{\gamma-1}$ функция $\psi_{\alpha_0}$ монотонна
на интервале $(\omega^-_{\alpha_0},\omega^+_{\alpha_0})$. Отсюда и из обусловленных
тождеством \eqref{eq:fa1} неравенств $f_{\alpha_0}'(\omega^\pm_{\alpha_0})\neq 0$
вытекает также равномерная положительность функции $\psi_{\alpha_0}$
на указанном интервале.

Дифференцируя с учетом последнего факта функции $I_\varepsilon$ при $\varepsilon>0$
на малой окрестности $U\ni\alpha_0$, устанавливаем справедливость соотношений
\begin{align*}
	I'_\varepsilon(\alpha)&\equiv-\int_{\omega^-_\alpha+\varepsilon}^%
		{\omega^+_\alpha-\varepsilon}\dfrac{t^{2\gamma}\,dt}{%
		2f_\alpha^{3/2}(t)}-\dfrac{(\omega^+_\alpha)^{2\gamma}}{%
		f_\alpha'(\omega^+_\alpha)\sqrt{f_\alpha(\omega^+_\alpha-
		\varepsilon)}}+\dfrac{(\omega^-_\alpha)^{2\gamma}}{%
		f_\alpha'(\omega^-_\alpha)\sqrt{f_\alpha(\omega^-_\alpha+
		\varepsilon)}}\\
	&=-\int_{\omega^-_\alpha+\varepsilon}^%
		{\omega^+_\alpha-\varepsilon}\dfrac{t^{2\gamma}\,dt}{%
		2f_\alpha^{3/2}(t)}-\left.\dfrac{A_\alpha(t)}{\psi_\alpha(t)
		\sqrt{f_\alpha(t)}}\right|^{\omega^+_\alpha-\varepsilon}_%
		{\omega^-_\alpha+\varepsilon}+O(\varepsilon^{1/2})\\
	&=-\int_{\omega^-_\alpha+\varepsilon}^{\omega^+_\alpha-\varepsilon}
		\left[\dfrac{t^{2\gamma}}{2f^{3/2}_\alpha(t)}+
		\left(\dfrac{A_\alpha(t)}{\psi_\alpha(t)\sqrt{f_\alpha(t)}}
		\right)'\right]\,dt+O(\varepsilon^{1/2}),
\end{align*}
где асимптотические оценки остаточных слагаемых равномерны по окрестности $U$.
Соответственно, функция $I_0$ также дифференцируема на $U$ и имеет производную
\begin{align*}
	I_0'(\alpha)&\equiv-\int_{\omega^-_\alpha}^{\omega^+_\alpha}\left[
		\dfrac{t^{2\gamma}}{2f^{3/2}_\alpha(t)}+\left(\dfrac{A_\alpha(t)}{%
		\psi_\alpha(t)\sqrt{f_\alpha(t)}}\right)'\right]\,dt\\
	&=-\int_{\omega^-_\alpha}^{\omega^+_\alpha}\dfrac{1}{\psi_\alpha^2(t)
		\sqrt{f_\alpha(t)}}\cdot\left[\dfrac{t^{2\gamma}\psi_\alpha^2(t)}{%
		2f_\alpha(t)}+{}\right.\\
	&\kern 2cm {}+\left(A'_\alpha(t)-\dfrac{A_\alpha(t)f'_\alpha(t)}{%
		2f_\alpha(t)}\right)\cdot\psi_\alpha(t)-A_\alpha(t)\psi'_\alpha(t)
		\biggr]\,dt\\
	&=-\int_{\omega^-_\alpha}^{\omega^+_\alpha}\dfrac{1}{\psi_\alpha^2(t)
		\sqrt{f_\alpha(t)}}\cdot\left[\dfrac{t^{2\gamma}\psi_\alpha^2(t)}{%
		2f_\alpha(t)}-{}\right.\\
	&\kern 2cm {}-\left(\dfrac{t^{2\gamma}\psi_\alpha(t)}{2f_\alpha(t)}+
		4\gamma\,(2\gamma-1)t^{2\gamma-2}f_\alpha(t)\right)\cdot
		\psi_\alpha(t)-A_\alpha(t)\psi'_\alpha(t)\biggr]\,dt\\
	&=\int_{\omega^-_\alpha}^{\omega^+_\alpha}\dfrac{1}{\psi_\alpha^2(t)
		\sqrt{f_\alpha(t)}}\cdot\biggl[4\gamma\,(2\gamma-1)
		t^{2\gamma-2}f_\alpha(t)\cdot[(f'_\alpha(t))^2-2f_\alpha(t)
		f_\alpha''(t)]-{}\\
	&\kern 2cm {}-2f_\alpha(t)f_\alpha'''(t)\cdot[t^{2\gamma}
		f_\alpha'(t)-4\gamma t^{2\gamma-1}f_\alpha(t)]\biggr]\,dt\\
	&=\int_{\omega^-_\alpha}^{\omega^+_\alpha}\dfrac{8\gamma\,(2\gamma-1)
		t^{2\gamma-3}f_\alpha(t)f'_\alpha(t)}{\psi_\alpha^2(t)
		\sqrt{f_\alpha(t)}}\,dt\\
	&=\int_{\omega^-_\alpha}^{\tau_\alpha}\dfrac{16\gamma\,(2\gamma-1)
		t^{2\gamma-3}}{3\psi_\alpha^2(t)}\,d(f_\alpha^{3/2})+
		\int_{\tau_\alpha}^{\omega^+_\alpha}\dfrac{16\gamma\,(2\gamma-1)
		t^{2\gamma-3}}{3\psi_\alpha^2(t)}\,d(f_\alpha^{3/2}).
\end{align*}
С учетом обусловленного \eqref{eq:psi} и \eqref{eq:fa3} убывания функции $t\mapsto
t^{2\gamma-3}\psi_{\alpha_0}^{-2}(t)$, это гарантирует справедливость неравенства
$I_0'(\alpha_0)>0$. Произвольность выбора точки $\alpha_0>\gamma^{-\gamma}\,
(1-\gamma)^{\gamma-1}$ означает потому строгое возрастание функции $I_0$ на всей
ее области определения. С другой стороны, при $\alpha\to\gamma^{-\gamma}\,
(1-\gamma)^{\gamma-1}+0$ имеют место асимптотические соотношения
\begin{align*}
	I_0(\alpha)&=\int_{\omega^-_\alpha}^{\omega^+_\alpha}
		\dfrac{d(f_\alpha')}{f_\alpha''\,\sqrt{f_\alpha}}\\
	&=\dfrac{1+o(1)}{2\sqrt{1-\gamma}}\cdot\left[\int_{\tau_\alpha}^%
		{\omega^-_\alpha}\dfrac{\sqrt 2\,d(f_\alpha')}{
		\sqrt{[f_\alpha'(\omega^-_\alpha)]^2-[f_\alpha']^2}}+
		\int_{\omega^+_\alpha}^{\tau_\alpha}\dfrac{\sqrt 2\,
		d(f_\alpha')}{\sqrt{[f_\alpha'(\omega^+_\alpha)]^2-
		[f_\alpha']^2}}\right]
\end{align*}
\begin{align*}
	&=\dfrac{1+o(1)}{\sqrt{1-\gamma}}\cdot\int_0^1\dfrac{\sqrt 2\,dt}{%
		\sqrt{1-t^2}}\\
	&=\dfrac{\pi}{\sqrt{2\cdot(1-\gamma)}}+o(1)\\ &=(1-\gamma)^{-1}+o(1).
\end{align*}
Тем самым, доказываемое утверждение справедливо.
\end{proof}

Доказательство предложения \ref{theorem2} завершено. В заключение отметим,
что приведенное доказательство предложения \ref{prop:pi2shao} основано
на идее доказательства основной леммы работы~\cite{Naz}.

\end{document}